\documentclass[11pt]{amsart}

\usepackage[frame,cmtip,arrow,matrix,line,graph,curve]{xy}
\usepackage{amsmath,amssymb, amscd}
\usepackage{amsfonts,amsxtra,amsthm}
\usepackage{mathrsfs}
\usepackage{eucal}
\usepackage{latexsym}

\numberwithin{equation}{section}

\newtheorem{theorem}{Theorem}[section]
\newtheorem{proposition}[theorem]{Proposition}
\newtheorem{corollary}[theorem]{Corollary}
\newtheorem{lemma}[theorem]{Lemma}

\newtheorem{remit}[theorem]{Remark}
\newtheorem{definit}[theorem]{Definition}

\newenvironment{definition}{\begin{definit}\rm}{\end{definit}}

\newcommand{\pp}{\mathbb{P}}
\newcommand{\qq}{\mathbb{Q}}
\newcommand{\cc}{\mathbb{C}}

\newcommand{\zz}{\mathbb{Z}}

\newcommand{\Hom}{\mathrm{Hom}}

\newcommand{\Gr}{\mathrm{Gr}}
\newcommand{\Ext}{\mathrm{Ext}}

\newcommand{\codim}{\mathrm{codim \,}}

\newcommand{\git}{/\!\!/}

\newcommand{\cE}{\mathcal{E} }
\newcommand{\cH}{\mathcal{H} }

\newcommand{\cM}{\mathcal{M} }
\newcommand{\cN}{\mathcal{N} }
\newcommand{\cL}{\mathcal{L} }
\newcommand{\cO}{\mathcal{O} }

\newcommand{\cV}{\mathcal{V} }
\newcommand{\cW}{\mathcal{W} }

\newcommand{\End}{\mathrm{End} }
\newcommand{\Pic}{\mathrm{Pic} }

\newcommand{\tcL}{\tilde{\mathcal{L}} }

\newcommand{\bM}{\mathbf{M} }
\newcommand{\bS}{\mathbf{S} }
\newcommand{\bR}{\mathbf{R} }

\begin{document}

\title[Stringy E-function of moduli of Higgs bundles]
{The stringy E-function of the moduli space of Higgs bundles with
trivial determinant }
\date{June, 2005}

\author{Young-Hoon Kiem and Sang-Bum Yoo}
\address{Department of Mathematics and Research Institute
of Mathematics, Seoul National University, Seoul 151-747, Korea}
\email{kiem@math.snu.ac.kr} \email{gauss76@snu.ac.kr}
\thanks{Partially supported by KRF}

\begin{abstract}
Let $\bM$ be the moduli space of semistable rank 2 Higgs pairs
$(V,\phi)$ with trivial determinant over a smooth projective curve
$X$ of genus $g\ge 2$. We compute the stringy E-function of $\bM$
and prove that there does not exist a symplectic desingularization
of $\bM$ for $g\ge 3$.
\end{abstract}
\maketitle


\section{Introduction}

Let $X$ be a smooth projective curve of genus $g\ge 2$. Let $\bM$
be the moduli space of semistable Higgs pairs $(V,\phi)$ over $X$
with $V$ a rank 2 vector bundle with $\det V\cong \cO_X$ and
$\phi\in H^0(\End_0V\otimes K_X)$. Then $\bM$ is a \emph{singular}
quasi-projective irreducible normal variety of dimension $6g-6$.
The locus $\bM^s$ of stable pairs in $\bM$ is an open dense subset
which is equipped with a symplectic form\footnote{In this paper, a
\emph{symplectic form} is always a \emph{holomorphic} 2-form which
is nondegenerate everywhere.} \cite{Hit1, Nit, simp}. In this
paper we give an explicit formula for the stringy E-function
$E_{st}(\bM)$ of $\bM$ as defined in \cite{Bat98} which retains
useful information about the singularities (Theorem 5.2).

The moduli space $\bM$ of Higgs pairs can be also thought of as
the moduli space of sheaves of certain topological type on the
\emph{symplectic} surface $T^*X$ \cite[\S6]{simp}. If a surface is
equipped with a symplectic form, it induces a natural symplectic
form on the smooth part of a moduli space of sheaves on the
surface \cite{Muk84}. A natural question raised by O'Grady
\cite{ogrady} asks whether there exists a desingularization of
such a moduli space on which the symplectic form extends
everywhere without degeneration. It was shown in \cite{KL04, KLS,
kiem04, CK1, CK2} that except for the two 10-dimensional moduli
spaces studied by O'Grady \cite{ogrady, ogrady2}, there does not
exist a symplectic desingularization when the surface is K3 or
Abelian.

In this paper we study the Kirwan desingularization of $\bM$ by
using O'Grady's analysis of the K3 surface case and show the
following.

\begin{theorem} (Theorem 4.1 and Corollary 5.4.) \\
(1) When $g=2$, there is a symplectic desingularization of $\bM$.
\\ (2) When $g\ge 3$, there does not exist a symplectic
desingularization of $\bM$.
\end{theorem}
If there exists a symplectic desingularization of $\bM$, then it
has to be a crepant resolution as the canonical bundle $K_\bM$ is
trivial. Hence the stringy E-function of $\bM$ has to be equal to
the Hodge-Deligne polynomial of the symplectic desingularization
which is a polynomial with integer coefficients (Theorem 2.1). For
the non-existence result (2), it suffices to prove that the
stringy E-function $E_{st}(\bM)$ of $\bM$ is not a polynomial for
$g\ge 3$. Because we have an explicit formula of $E_{st}(\bM)$
(Theorem 5.2), it is a simple matter to prove this.

To compute $E_{st}(\bM)$ we consider all possible types of
semistable pairs, and find a description of the locally closed
subvariety of $\bM$ corresponding to each type. Then we compute
the Hodge-Deligne polynomials of the subvarieties explicitly.

In \S2, we recall basic facts about stringy E-function and Higgs
pairs. In \S3, we compute the Hodge-Deligne polynomial of the
stable locus $\bM^s$. In \S4, we construct the Kirwan
desingularization of $\bM$ along the line of O'Grady's analysis in
\cite{ogrady} and show that $\bM$ admits a symplectic
desingularization when $g=2$. In \S5, we complete the computation
of $E_{st}(\bM)$ and prove non-existence of symplectic
desingularization.

In the context of Mirror Symmetry, Hausel and
Thaddeus computed the stringy E-function of the moduli space of
Higgs pairs of odd degree \cite{HaTa}. All the varieties in this
paper are defined over the complex number field.

It is our pleasure to express gratitude to Professors S. Ramanan,
N. Niture, C. Sorger and J.-M. Hwang for useful discussions.


\section{Preliminaries}

In this section we collect some facts that we shall use in this
paper.

\vspace{.5cm}\textbf{2.1. Stringy E-function}
\bigskip

Stringy E-function introduced in \cite{Bat98} is a new invariant
of varieties which retains useful information about singularities.
We recall the definition and basic facts about stringy E-functions
from \cite{Bat98,DL99}. Let $W$ be a normal irreducible variety
with at worst log-terminal singularities, i.e.
\begin{enumerate} \item W is $\qq$-Gorenstein;
\item for a resolution of singularities $\rho: V\to W$ such that
the exceptional locus of $\rho$ is a divisor $D$ whose irreducible
components $D_1,\cdots,D_r$ are smooth divisors with only normal
crossings, we have \[K_V=\rho^*K_W+\sum^r_{i=1} a_iD_i \] with
$a_i>-1$ for all $i$, where $D_i$ runs over all irreducible
components of $D$. The divisor $\sum^r_{i=1}a_iD_i$ is called the
\textit{discrepancy divisor}. \end{enumerate}

For each subset $J\subset I=\{1,2,\cdots,r\}$, define
$D_J=\cap_{j\in J}D_j$, $D_\emptyset =V$ and $D^0_J=D_J-\cup_{i\in
I-J}D_i$. Then the stringy E-function of $W$ is defined by
\begin{equation} \label{eqn:stringy E-function}
E_{st}(W;u,v)=\sum_{J\subset I}E(D^0_J;u,v)\prod_{j\in
J}\frac{uv-1}{(uv)^{a_j+1}-1} \end{equation} where \[ E(Z;u,v) =
\sum_{p,q}\sum_{k\geq 0} (-1)^kh^{p,q}(H^k_c(Z;\cc))u^pv^q \] is
the Hodge-Deligne polynomial for a variety $Z$. We will also use
the alias, \emph{E-polynomial}, for the Hodge-Deligne polynomial
and often use the abbreviation $E(Z)$ for $E(Z;u,v)$. Note that
the Hodge-Deligne polynomials have
\begin{enumerate}
\item the additive property: $E(Z;u,v)=E(U;u,v)+E(Z-U;u,v)$ if $U$
is a smooth open subvariety of $Z$; \item the multiplicative
property: $E(Z;u,v)=E(B;u,v)E(F;u,v)$ if $Z$ is a Zariski locally
trivial $F$-bundle over $B$. \end{enumerate}

By \cite[Theorem 6.27]{Bat98}, the function $E_{st}$ is
independent of the choice of a resolution (Theorem 3.4 in
\cite{Bat98}) and the following holds.
\begin{theorem} \label{thm:Batyrev's result} (\cite[Theorem
3.12]{Bat98}) Suppose $W$ is a $\qq$-Gorenstein algebraic variety
with at worst log-terminal singularities. If $\rho:V\to W$ is a
crepant desingularization (i.e. $\rho^*K_W=K_V$) then
$E_{st}(W;u,v)=E(V;u,v)$. In particular, $E_{st}(W;u,v)$ is a
polynomial.
\end{theorem}

\vspace{.5cm}\textbf{2.2. Higgs pairs}
\bigskip

Let $X$ be a smooth projective curve of genus $g\ge 2$. A Higgs
pair, more precisely an $SL(2)$-Higgs pair, is a pair of a rank 2
vector bundle $V$ with trivial determinant and a section $\phi$ of
$\End_0V\otimes K_X$ where $K_X$ is the canonical bundle of $X$
and $\End_0V$ denotes the traceless part of $\End V$. To construct
the moduli space of Higgs pairs, a stability condition has to
imposed and Hitchin introduced the following.
\begin{definition}
(1) A Higgs pair $(V,\phi)$ with $\det V\cong \cO_X$ is stable
(resp. semistable) if for any nonzero proper subbundle $W$
satisfying $\phi (W)\subset W\otimes K$, we have $\deg W<0$ (resp.
$\deg W\le 0$).

(2) A Higgs pair $(V,\phi)$ with $\det V\cong \cO_X$ is polystable
if it is either stable or a direct sum $(L,\psi)\oplus
(L^{-1},-\psi)$ where $L\in \Pic^0(X)$ is a line bundle of degree
0 and $\psi$ is a section of $\Hom(L,L)\otimes K_X\cong K_X$.
\end{definition}

If a bundle $V$ is stable (resp. semistable),\footnote{ A vector
bundle $V$ with trivial determinant is stable (resp. semistable)
for any nonzero proper subbundle $W$, we have $\deg W<0$ (resp.
$\deg W\le 0$). A vector bundle is strictly semistable if it is
semistable but not stable.}
 the Higgs
pair $(V,\phi)$ is stable (resp. semistable) for any choice of
$\phi\in H^0(\End_0V\otimes K_X)$.

The set of isomorphisms classes of polystable pairs $(V,\phi)$
admits a structure of quasi-projective variety of dimension $6g-6$
\cite{Hit1, Nit, simp} and we denote it by $\bM$. Furthermore, it
is known \cite{Hit1, simp} that $\bM$ is an irreducible normal
variety. The locus $\bM^s$ of stable pairs $(V,\phi)$ is a smooth
open dense subvariety and its complement is precisely the locus of
singularities, isomorphic to $T^*J/\zz_2$ where $J=\Pic^0(X)$ is
the Jacobian and $\zz_2$ acts on the cotagent bundle $T^*J=J\times
H^0(K_X)$ by $(L,\psi)\mapsto (L^{-1},-\psi)$.

The stable locus $\bM^s$ is homeomorphic to the space of
irreducible representations of the fundamental group of $X$ into
$SL(2)$ and the complex structures of $X$ and $SL(2)$ induce two
integrable complex structures on $\bM^s$. Therefore, $\bM^s$
admits a hyperk\"ahler metric and thus there is a symplectic form
on $\bM^s$. In particular the canonical bundle of $\bM$ is trivial
and $\bM$ is Gorenstein.\footnote{A normal variety $\bM$ is
Gorenstein if the canonical divisor $K_{\bM}$ is Cartier
\cite{Bat98}. In our case, the existence of symplectic form
guarantees $K_{\bM^s}=0$ and hence $K_\bM=0$ since $\codim
(\bM-\bM^s)\ge 2$. Obviously $0$ is Cartier.} We will see in
Theorem 4.1 that $\bM$ has only log terminal singularities.
Therefore, the stringy E-function of $\bM$ is a well-defined
rational function which we intend to compute.


\section{Stable pairs}

In this section we compute the E-polynomial of the stable part
$\bM^s$. Let $(V,\phi)\in \bM^s$. There are three possibilities
for $V$ : \begin{enumerate}\item $V$ is stable
\item $V$ is strictly semistable
\item $V$ is unstable (= not semistable)\end{enumerate} We deal
with these cases separately in the subsequent subsections.

\vspace{.5cm}\textbf{3.1. Stable case}
\bigskip

Let $\cN^s$ be the moduli space of rank 2 stable bundles with
trivial determinant over $X$. $\cN^s$ is a $3g-3$ dimensional
quasi-projective variety and the Hodge-Deligne polynomial of
$\cN^s$ is
\begin{equation}\label{eq3.0}
 E(\cN^s)=\frac{(1-u^2v)^g(1-uv^2)^g-(uv)^{g+1}(1-u)^g(1-v)^g}{(1-uv)(1-(uv)^2)}\end{equation}
$$-\frac12(\frac{(1-u)^g(1-v)^g}{1-uv}+\frac{(1+u)^g(1+v)^g}{1+uv})
$$
from \cite[(18)]{kiem} or \cite[\S6.2]{KL}.

If $V$ is stable, a pair $(V,\phi)$ is stable for any $\phi\in
H^0(\End_0V\otimes K_X)$. It is well-known that such pairs are
parameterized by the cotangent bundle $T^*\cN^s$ which is embedded
in $\bM^s$ as an open subvariety \cite{Hit1}. Hence the
E-polynomial of the locus of stable pairs $(V,\phi)$ with $V$
stable is
\begin{equation}\label{eq3.1}
E(T^*\cN^s)=(uv)^{3g-3}\big(
\frac{(1-u^2v)^g(1-uv^2)^g-(uv)^{g+1}(1-u)^g(1-v)^g}{(1-uv)(1-(uv)^2)}\end{equation}
$$-\frac12(\frac{(1-u)^g(1-v)^g}{1-uv}+\frac{(1+u)^g(1+v)^g}{1+uv}) \big)
$$

\vspace{.5cm} \textbf{3.2. Strictly semistable case}
\bigskip

When $V$ is a strictly semistable rank 2 bundle with trivial
determinant, there are four possibilities:
\begin{enumerate}
\item[Type I:] $V=L\oplus L^{-1}$ for $L\in \Pic^0(X)$ with
$L\ncong L^{-1}$
\item[Type II:] $V$ is a nontrivial extension of $L^{-1}$ by $L$
for $L\in \Pic^0(X)$ with $L\ncong L^{-1}$
\item[Type III:] $V=L\oplus L^{-1}$ for $L\in \Pic^0(X)$ with
$L\cong L^{-1}$
\item[Type IV:] $V$ is a nontrivial extension of $L^{-1}$ by $L$
for $L\in \Pic^0(X)$ with $L\cong L^{-1}$
\end{enumerate}

We consider the loci of stable pairs $(V,\phi)$ for the above four
cases separately.

\vspace{.5cm} \textit{3.2.1. Type I}
\bigskip

Let $J=\Pic^0(X)$ and $J_0$ be the locus of $L\in J$ with $L\cong
L^{-1}$ so that $J_0\cong \zz_2^{2g}$. Let $J^0=J-J_0$ be the
complement of $J_0$ in $J$. Then the bundle $V$ is parameterized
by $J^0/\zz_2=J/\zz_2-J_0$ where $-1\in \zz_2$ acts as $L\mapsto
L^{-1}$.

When $V=L\oplus L^{-1}$ for $L\in J^0$, we have the decomposition
\[
H^0(\End_0(V)\otimes K_X)=H^0(K_X)\oplus H^0(L^2K_X)\oplus
H^0(L^{-2}K_X).
\]
For $\phi\in H^0(\End_0(V)\otimes K_X)$, write $\phi=(a,b,c)$  or
\[\phi=\left( \begin{matrix} a&b\\ c& -a\end{matrix}\right)\]
with $a\in H^0(K_X)$, $b\in H^0(L^2K_X)$ and $c\in
H^0(L^{-2}K_X)$. Then a pair $(V,\phi)$ is stable if and only if
$L$ and $L^{-1}$ are not preserved by $\phi$, i.e. $b\ne 0$ and
$c\ne 0$.

Because the automorphism group of $V=L\oplus L^{-1}$ with $L\ncong
L^{-1}$ is $\cc^*\times \cc^*$, $(V,\phi_1)\cong (V,\phi_2)$ for
$\phi_i=(a_i,b_i,c_i)$ if and only if
\[
\phi_1=\left( \begin{matrix} t&0\\ 0&t^{-1}\end{matrix}\right)
\phi_2 \left( \begin{matrix} t^{-1}&0\\
0&t\end{matrix}\right)
\] for some $t\in \cc^*$, i.e. $a_1=a_2$, $b_1=t^2b_2$ and $c_1=t^{-2}c_2$.
Therefore, for fixed $V=L\oplus L^{-1}$ of Type I, the isomorphism
classes of stable pairs $(V,\phi)$ are parameterized by
\begin{equation}\label{eq3.2}
\begin{array}{ll}&H^0(K_X)\times \frac{(H^0(L^2K_X)-0)\times
(H^0(L^{-2}K_X)-0)}{\cc^*}\\ &=\cc^g\times
\frac{(\cc^{g-1}-0)\times (\cc^{g-1}-0)
}{\cc^*}.\end{array}\end{equation} It is well-known that the
quotient $\pp (\cc^{g-1}\times \cc^{g-1})\git \cc^*$ is
$\pp^{g-2}\times \pp^{g-2}$ and $\cc^{g-1}\times\cc^{g-1}\git
\cc^*$ is the affine cone over $\pp^{g-2}\times \pp^{g-2}$ while
$(0\times \cc^{g-1})\cup (\cc^{g-1}\times 0)$ is the inverse image
in $\cc^{g-1}\times \cc^{g-1}$ of the vertex of the affine cone
via the quotient map. Hence
\[
\frac{(\cc^{g-1}-0)\times (\cc^{g-1}-0) }{\cc^*}\] is the line
bundle $\cO_{\pp^{g-2}\times\pp^{g-2}}(-1,-1)$ minus the zero
section and the E-polynomial of \eqref{eq3.2} is
\begin{equation}\label{eq3.3}
E(\cc^g)\cdot E(\cc^*)\cdot E(\pp^{g-2}\times \pp^{g-2})
=(uv)^g(uv-1)\left(
\frac{(uv)^{g-1}-1}{uv-1}\right)^2.\end{equation}

Now we let $V$ vary. Let $\cL\to J^0\times X$ be a universal
bundle of degree $0$ line bundles and let $\pi_1$, $\pi_2$ be the
projections of $J^0\times X$ to $J^0$ and $X$ respectively. Let
\[
\cW_j={\pi_1}_*(\cL^{2j}\otimes \pi_2^*K_X)
\]
for $j=0,1,-1$ and let $\cW=\cW_0\oplus \cW_1\oplus \cW_{-1}$.
Since $H^1(L^{\pm 2}K_X)\cong H^0(L^{\mp 2})^\vee =0$, $\cW_0,
\cW_{\pm 1}$ are vector bundles over $J^0$ of rank $g$ and $g-1$
respectively whose fibers over $L\in J^0$ are $H^0(K_X)$ and
$H^0(L^{\pm 2}K_X)$.

There is an obvious family of Higgs pairs $(\cV,\Phi)$
parameterized by $\cW$ with $\cV=\tcL\oplus \tcL^{-1}$ where
$\tcL$ is the pull-back of $\cL$ to $\cW\times X$ and this family
restricted to
\[ \cW_0\oplus (\cW_1-0)\oplus (\cW_{-1}-0)\]
parameterizes stable pairs where $0$ denotes the zero section.
Hence we have a morphism
\begin{equation}\label{eq3.5}
\cW_0\oplus (\cW_1-0)\oplus (\cW_{-1}-0)\to \bM^s
\end{equation}
 Since the action of $\cc^*$ on $\cW_j$
with weight $2j$ for $j=0,1,-1$ preserves the isomorphism classes
of stable pairs, \eqref{eq3.5} factors through
\begin{equation}\label{eq3.6}
\cW_0\oplus \frac{(\cW_1-0)\oplus (\cW_{-1}-0)}{\cc^*}\to \bM^s.
\end{equation}
Furthermore, the $\zz_2$-action on $J^0$ which interchanges $L$
with $L^{-1}$ obviously extends to $\cW$ interchanging $\cW_1$ and
$\cW_{-1}$. Thus we get a morphism
\begin{equation}\label{eq3.7}
\left[\cW_0\oplus \frac{(\cW_1-0)\oplus
(\cW_{-1}-0)}{\cc^*}\right]/\zz_2\to \bM^s.
\end{equation}
From construction, it is clear that this is a bijection onto the
locus of stable pairs $(V,\phi)$ with $V$ of Type I.

\begin{lemma}\label{lem1} The morphism \eqref{eq3.7} is an
isomorphism onto a locally closed subvariety of $\bM^s$.
\end{lemma}
\begin{proof}
First observe that the locus of $(V,\phi)$ with $V$ of Type I is
locally closed. Indeed, given a family of stable pairs $(V,\phi)$
on $X$ parameterized by a variety $T$, the locus in $T$ of
semistable $V$ is open and the locus of decomposable semistable
$V$ is closed in this open set. The condition $L\ncong L^{-1}$ is
determines an open subset of this locally closed set.

Let $(\tilde{\cV},\tilde{\Phi})$ be a family of stable pairs
$(V,\phi)$ with $V$ of Type I parameterized by a variety $T$. Then
there is a morphism $T\to J^0/\zz_2$ which sends $V=L\oplus
L^{-1}$ to the $\zz_2$-orbit $(L,L^{-1})$. Indeed, as
$\tilde{\cV}$ is a family of semistable bundles, there is a
morphism $T\to \cN$ whose image is obviously the Kummer variety
$J/\zz_2$ minus $J_0$, i.e. $J^0/\zz_2$.

Moreover, the section $\tilde{\Phi}$ of $\End_0(\tilde{\cV}\otimes
K_X)$ induces a lifting of the morphism $T\to J^0/\zz_2$ to a
morphism
\[
T\to \left[\cW_0\oplus \frac{(\cW_1-0)\oplus
(\cW_{-1}-0)}{\cc^*}\right]/\zz_2.
\]
Indeed, if we let $\hat{T}$ be the fiber product of $T$ and $J^0$
over $J^0/\zz_2$ and ${f}:\hat{T}\to J^0$ be obvious map, then the
pull-back $\hat{\cV}$ of $\tilde{\cV}$ to $\hat{T}\times X$ is
$(f\times 1_X)^*(\cL\oplus \cL^{-1})$ and from the commutative
square
\[\xymatrix{
\hat{T}\times X \ar[r]^{f\times 1}\ar[d]_{\pi}& J^0\times X\ar[d]^{\pi'}\\
\hat{T}\ar[r]^f &J^0 }\] where the vertical maps are the
projections, $\tilde{\Phi}$ induces a section of
$$\pi_*(\End_0(\hat{\cV})\otimes K_X)=\pi_*(f\times
1_X)^*(\End_0(\cL\oplus \cL^{-1})\otimes K_X)$$
$$=f^*\pi'_*\End_0(\cL\oplus \cL^{-1})\otimes K_X=f^*\cW.$$
Therefore, $\tilde{\Phi}$ gives us a morphism
$$\hat{T}\to \cW_0\oplus (\cW_1-0)\oplus (\cW_{-1}-0)$$
as $\hat{\cV}$ is a family of stable pairs. So we obtain a
morphism
\[
T\to \left[\cW_0\oplus \frac{(\cW_1-0)\oplus
(\cW_{-1}-0)}{\cc^*}\right]/\zz_2
\]
after taking quotients. Obviously this gives us the inverse of
\eqref{eq3.7}.
\end{proof}

Consequently, the locus of stable pairs $(V,\phi)$ with $V$ of
Type I is
\[
\left[\cW_0\oplus \frac{(\cW_1-0)\oplus
(\cW_{-1}-0)}{\cc^*}\right]/\zz_2
\]
which is a fiber bundle over $J^0/\zz_2$ with fiber \eqref{eq3.2}.
The E-polynomial of the fiber bundle is the E-polynomial of the
$\zz_2$-invariant part of
$$H^*_c\left(\cW_0\oplus \frac{(\cW_1-0)\oplus
(\cW_{-1}-0)}{\cc^*}\right).$$ From the Leray spectral sequence,
this is precisely
\begin{equation}\label{eq3.91}
(uv)^g\cdot (uv-1)\cdot \left[ E(\pp^{g-2}\times \pp^{g-2})^+\cdot
E(J^0)^+ + E(\pp^{g-2}\times \pp^{g-2})^-\cdot E(J^0)^- \right]
\end{equation}
where $E(Z)^\pm$ denote the E-polynomials of the $\zz_2$-invariant
and anti-invariant part of the compact support cohomology of a
variety $Z$. From the computation of $E(\tilde{D}_2^{(2)})$ in
\cite[p516]{KL}, we deduce that \eqref{eq3.91} is
\begin{equation}\label{eq3.10}
(uv)^g\cdot (uv-1)\cdot \big[ \big(\frac12 (1-u)^g(1-v)^g+\frac12
(1+u)^g(1+v)^g
-2^{2g}\big)\frac{((uv)^g-1)((uv)^{g-1}-1)}{(uv-1)((uv)^2-1)}\end{equation}
$$+\big(\frac12 (1-u)^g(1-v)^g-\frac12 (1+u)^g(1+v)^g
\big)uv\frac{((uv)^{g-1}-1)((uv)^{g-2}-1)}{(uv-1)((uv)^2-1)} \big]
.$$ This is the E-polynomial of the locus of stable pairs with $V$
of Type I.

\vspace{.5cm} \textit{3.2.2. Type II}
\bigskip

Now we consider the locus of stable pairs $(V,\phi)$ with $V$ of
Type II. Let $\cL\to J^0\times X$ be a universal line bundle on
$J^0$ and let $\pi:J^0\times X\to J^0$ be the projection. Then
$R^1\pi_*\cL^2$ is a vector bundle of rank $g-1$ by Riemann-Roch.
Let $\Lambda=\pp R^1\pi_*\cL^2$ be the projectivization of
$R^1\pi_*\cL^2$ and let $\cL^\#$ be the pull-back of $\cL$ by
$\Lambda\times X\to J^0\times X$. It is well-known that there is a
universal extension bundle
\begin{equation}\label{eq3.11}
0\to \cL^\#\otimes \cO_\Lambda(1)\to \cV^\# \to (\cL^\#)^{-1}\to
0\end{equation} over $\Lambda\times X$ where $\cO_\Lambda (1) $ is
the hyperplane bundle over the projective bundle $\Lambda$. Then
the family $\cV^\#$ over $\Lambda\times X$ parameterizes all
isomorphism classes of rank 2 vector bundles $V$ of Type II.

Next we consider the bundle $\cE nd_0\cV^\#\otimes K_X $ where
$K_X$ denotes the pull-back of the canonical bundle of $X$ by
abuse of notation. This fits into a short exact sequence
\cite[(3.7)]{Hit1}
\begin{equation}\label{eq3.12}
0\to \cH om(\cV^\#, K_X\cL^\#)\to \cE nd_0\cV^\#\otimes K_X\to
K_X(\cL^\#)^{-2}\to 0
\end{equation}
which gives rise to an exact sequence
\begin{equation}\label{eq3.13}
0\to p_*\cH om(\cV^\#, K_X\cL^\#)\to p_*\cE nd_0\cV^\#\otimes
K_X\to p_*K_X(\cL^\#)^{-2}\to R^1p_*\cH om(\cV^\#, K_X\cL^\#)
\end{equation}
where $p:\Lambda\times X\to \Lambda$ is the projection. From
\eqref{eq3.11} we also have an exact sequence
\[
R^1p_*K_X(\cL^\#)^2\to R^1p_*\cH om(\cV^\#,K_X\cL^\#)\to
R^1p_*K_X\cO_\Lambda (-1) \to 0.
\]
By the Serre duality, $R^1p_*K_X(\cL^\#)^2=0$ and
$R^1p_*K_X\otimes \cO_\Lambda (-1)$ is a line bundle over
$\Lambda$ since $X$ is irreducible projective of dimension 1.
Hence the last map in \eqref{eq3.13} gives us
\[
p_*K_X(\cL^\#)^{-2}\to R^1p_*\cH om(\cV^\#, K_X\cL^\#)\cong
R^1p_*K_X\cO_\Lambda (-1).
\]
Over a point $s\in \Lambda$ lying over $[L]\in J^0$, the fibers
are
\[
H^0(K_XL^{-2})\to H^1(V^*\otimes K_XL)\cong H^1(K_X)\cong \cc
\]
which is the multiplication by the extension class of $V$, a
nonzero representative in $H^1(L^2)$ of $s$. In particular, this
is surjective and hence the kernel of the last map in
\eqref{eq3.13} is a vector bundle of rank $g-2$ because the rank
of $p_*K_X(\cL^\#)^{-2}$ is $g-1$.

From \eqref{eq3.11} again, we have an exact sequence
\begin{equation}\label{eq3.14}
0\to p_*K_X(\cL^\#)^2\to p_*\cH om(\cV^\#,K_X\cL^\#)\to
p_*K_X\cO_\Lambda (-1) \to R^1p_*K_X(\cL^\#)^2=0.
\end{equation}
The first term in \eqref{eq3.14} is a vector bundle over $\Lambda$
of rank $g-1$ and the last term is a vector bundle of rank $g$.
Therefore, $p_*(\cE nd_0\cV^\#\otimes K_X)$ is a vector bundle
over $\Lambda$ of rank $3g-3$ from \eqref{eq3.13}.

Let $(\cV^\dagger,\Phi^\dagger)$ be the obvious family of Higgs
pairs over $p_*(\cE nd_0\cV^\#\otimes K_X)\times X$ where
$\cV^\dagger$ is the pull-back of $\cV^\#$ via the bundle
projection $p_*(\cE nd_0\cV^\#\otimes K_X)\times X\to
\Lambda\times X$. A pair $(V,\phi)$ in this family lying over
$[L]\in J^0$ is stable if and only if $L$ is not preserved by
$\phi$. This amounts to saying that the image of this point by the
middle map in \eqref{eq3.13}
\[
p_*\cE nd_0\cV^\#\otimes K_X\to p_*K_X(\cL^\#)^{-2}
\]
is nonzero. Therefore the stable locus in $p_*(\cE
nd_0\cV^\#\otimes K_X)$ is precisely the complement of the
subbundle $p_*\cH om(\cV^\#, K_X\cL^\#)$. For a bundle $V$ in
$\cV^\#$, the automorphism group is trivial. So we obtain the
following.

\begin{lemma} The locus of stable pairs $(V,\phi)$ in
$\bM^s$ with $V$ of Type II is locally closed and isomorphic to
\begin{equation}\label{eq3.16}
p_*(\cE nd_0\cV^\#\otimes K_X)-p_*\cH om(\cV^\#, K_X\cL^\#).
\end{equation} In particular the E-polynomial of this locus is
\begin{equation}\label{eq3.17}
[(uv)^{3g-3}-(uv)^{2g-1}]\cdot
E(\Lambda)=[(uv)^{3g-3}-(uv)^{2g-1}] \frac{(uv)^{g-1}-1}{uv-1}
[(1-u)^g(1-v)^g-2^{2g}].
\end{equation}\end{lemma}

The proof is similar but much easier than Lemma \ref{lem1} and we
omit it.

\vspace{.5cm}\textit{3.2.3. Type III}
\bigskip

We turn now to the Type III case. Let $L\cong L^{-1}\in J_0\cong
\zz_2^{2g}$. By tensoring $L\in J_0$, we may restrict our concern
to the case $L=\cO_X$ so that $V=\cO_X\oplus \cO_X$.

Since $H^0(\End _0V\otimes K_X)=H^0(K_X)\otimes sl(2)\cong
\cc^g\otimes sl(2)$, a Higgs field is of the form
\[
\phi=\left(\begin{matrix} a&b\\ c&-a\end{matrix}\right)\qquad
\text{for }a,b,c\in H^0(K_X).
\]
If the pair $(V,\phi)$ is not stable, there must be an injective
map
\[
\imath:\cO_X\to \cO_X\oplus \cO_X=V
\]
and $\phi$ should preserve the image of $\imath$. Let $p$ (resp.
$q$) in $\cc$ be the composition of $\imath$ and the projection
$V\to \cO_X$ onto the first (resp. second) component. If $\phi$
preserves the image of $\imath$, there exists $\lambda\in
H^0(K_X)$ such that
\[ p\lambda=pa+qb,\quad q\lambda=pc-qa
\]
hold at the same time. It is an elementary exercise to show that
this condition is equivalent to saying that $\phi$ is conjugate to
an upper triangular matrix $\left(\begin{matrix} a'&b'\\
0&-a'\end{matrix}\right)$. By the Hilbert-Mumford criterion, we
conclude that $(V,\phi)$ is stable if and only if $\phi$ is
nonzero and the line $[\phi]\in \pp (\cc^g\otimes sl(2))$ is
stable with respect to the adjoint action of $SL(2)$. The
automorphism group of $V$ is $GL(2)$ which acts on $H^0(\End
_0V\otimes K_X)=\cc^g\otimes sl(2)$ by conjugation. Note that the
center $\cc^*$ of $GL(2)$ acts trivially. So we obtain the
following.

\begin{lemma} The locus of stable pairs $(V,\phi)$ in $\bM^s$ with
$V$ of Type III is the disjoint union of $2^{2g}$ locally closed
subvarieties, each isomorphic to
\[
(\cc^g\otimes sl(2))^{st}/SL(2)
\]
where $(\cc^g\otimes sl(2))^{st}$ is the set of nonzero $\phi\in
\cc^g\otimes sl(2)$ with $[\phi]\in \pp (\cc^g\otimes sl(2))$
stable. \end{lemma}

In particular, the E-polynomial of the locus of Type III is
\begin{equation}\label{eq3.171}
2^{2g} (uv-1)\cdot E(\pp(\cc^g\otimes sl(2))^{st}/SL(2)).
\end{equation}
We can compute $E(\pp(\cc^g\otimes sl(2))^{st}/SL(2))$ by Kirwan's
algorithm \cite{k2} as follows. (See \cite[\S4]{kiem}.) We start
with the $SL(2)$-equivariant cohomology of $\pp (\cc^g\otimes
sl(2))$ whose Hodge-Deligne series is
\begin{equation}\label{eq3.18}
\frac{1}{1-(uv)^2}\cdot \frac{1-(uv)^{3g}}{1-uv}
\end{equation}
and subtract out the Hodge-Deligne series of the unstable part
\begin{equation}\label{eq3.19}
(uv)^{2g-1} \cdot \frac{1+uv+\cdots +(uv)^{g-1}}{1-uv}.
\end{equation}
Next we blow up along $SL(2) \pp (\cc^g\otimes
\left(\begin{matrix}1&0\\ 0& -1\end{matrix}\right) )$ and delete
the unstable part. For the Hodge-Deligne series, we have to add
\begin{equation}\label{eq3.20}
\frac{1+uv+\cdots +(uv)^{g-1}}{1-(uv)^2} \cdot (uv+\cdots
+(uv)^{2g-3})\end{equation} $$ - (uv)^{g-1}\cdot \frac{1+uv+\cdots
+(uv)^{g-2}}{1-uv}\cdot (1+uv+\cdots +(uv)^{g-1}).$$ The
$SL(2)$-quotient of this blow-up is Kirwan's partial
desingularization whose E-polynomial\footnote{Note that this is
projective and the compact support cohomology is the same as the
ordinary cohomology which is isomorphic to the equivariant
cohomology of the stable part.} is
\begin{equation}\label{eq3.21}
\frac{(1-(uv)^{g-1})(1-(uv)^{g})(1-(uv)^{g+1})}{(1-uv)^2(1-(uv)^{2})}
\end{equation}
from \eqref{eq3.18}, \eqref{eq3.19}, and \eqref{eq3.20}. The
quotient of the exceptional divisor of the blow-up is a Zariski
locally trivial bundle over $\pp^{g-1}$ with fiber
$\pp^{g-2}\times_{\zz_2}\pp^{g-2}$ whose E-polynomial is
\begin{equation}\label{eq3.22}
\frac{1-(uv)^g}{1-uv}\cdot \frac12
\left[\big(\frac{1-(uv)^{g-1}}{1-uv}\big)^2+\frac{1-(uv)^{2g-2}}{1-(uv)^2}
\right]\end{equation} Upon subtracting \eqref{eq3.22} from
\eqref{eq3.21}, we deduce that
\[
E(\pp(\cc^g\otimes sl(2))^{st}/SL(2))=
\frac{(uv)^g(1-(uv)^{g-1})(1-(uv)^{g})}{(1-uv)(1-(uv)^{2})}.
\]
From \eqref{eq3.171}, the E-polynomial of the locus of Type III is
finally
\begin{equation}\label{eq3.23}
2^{2g}\cdot\frac{(uv)^g((uv)^{g-1}-1)((uv)^{g}-1)}{(uv)^2-1}.
\end{equation}

\vspace{.5cm} \textit{3.2.4. Type IV}
\bigskip

As in the Type III case, we may assume $L\cong \cO_X$ and $V$ is a
nontrivial extension of $\cO_X$ by $\cO_X$.  The isomorphism
classes of such bundles $V$ are parameterized by $$\Gamma=\pp
\Ext^1(\cO_X,\cO_X)=\pp H^1(\cO_X)\cong \pp^{g-1}.$$ There is a
universal extension bundle
\begin{equation}\label{eq3.24}
0\to \cO_\Gamma (1)\to \cV\to \cO_\Gamma \to 0
\end{equation}
where $\cO_\Gamma(1)$ is the hyperplane bundle over $\Gamma$. Let
$p:\Gamma\times X\to \Gamma$ be the projection. Exactly as in the
Type II case, we have an exact sequence of vector bundles
\begin{equation}\label{eq3.25}
0\to p_*\cH om(\cV, K_X)\to p_*(\cE nd_0\cV\otimes
K_X)\to\end{equation}
\[
\to p_*K_X\to R^1p_*\cH om (\cV,K_X)\cong R^1p_*K_X\otimes
\cO_\Gamma(-1)
\]
and the last map is surjective because the extensions $V$ are
nontrivial. Hence the kernel of the last map is a vector bundle
over $\Gamma$ of rank $g-1$. From \eqref{eq3.24}, we also have an
exact sequence of vector bundles
\begin{equation}\label{eq3.26}
0\to p_*K_X\to p_*\cH om(\cV,K_X)\to p_*K_X\otimes \cO_\Gamma
(-1)\to R^1p_*K_X\end{equation} whose last map is surjective
because the extensions $V$ are nontrivial. So the rank of the
vector bundle $p_*\cH om(\cV,K_X)$ over $\Gamma$ is $2g-1$.

As in the Type II case, if we consider the obvious family of Higgs
pairs over $p_*(\cE nd_0\cV\otimes K_X)\times X$, the locus of
stable pairs in $p_*(\cE nd_0\cV\otimes K_X) $ is precisely the
complement of the subbundle $p_*\cH om(\cV, K_X)$.

Finally the automorphisms of $V$ should be taken into account. For
a nontrivial extension $V$ of $\cO_X$ by $\cO_X$, the automorphism
group is the additive group $(\cc,+)$ and locally $q\in\cc$ acts
by \begin{equation}\label{eq3.25a}
q\cdot \phi=\left(\begin{matrix}a+qc & b-2qa-q^2c\\
c& -a-qc
\end{matrix}\right)\quad \text{ for }\phi=\left(\begin{matrix} a&b\\
c&-a\end{matrix}\right).\end{equation}

\begin{lemma}\label{lem3.4} The locus of stable pairs $(V,\phi)$ with $V$ of
Type IV is the disjoint union of $2^{2g}$ locally closed
subvarieties, each isomorphic to a $\cc^{g}$-bundle over a
$\cc^{g-2}$-bundle over a $(\cc^{g-1}-0)$-bundle over $\pp^{g-1}$.
All the bundles are Zariski locally trivial.\end{lemma}
\begin{proof}
Let $A$ be the kernel of the last map in \eqref{eq3.25} minus the
zero section. Then $A$ is a $(\cc^{g-1}-0)$-bundle over
$\Gamma\cong \pp^{g-1}$. We think of $$p_*(\cE nd_0\cV\otimes
K_X)-p_*\cH om(\cV,K_X)$$ as a vector bundle of rank $2g-1$ over
$A$.

The kernel of the last map in \eqref{eq3.26} gives rise to a
vector bundle $\mathfrak{A}$ over $A$ of rank $g-1$ and the second
map in \eqref{eq3.26} lifts to a $\cc$-equivariant map
\begin{equation}\label{eq3.25b}
\left[ p_*(\cE nd_0\cV\otimes K_X)-p_*\cH om(\cV,K_X)\right]\to
\mathfrak{A}
\end{equation}
of vector bundles over $A$ whose kernel is of rank $g$. The action
of $\cc$ on $\mathfrak{A}$ is linear $a\mapsto a+qc$ as the first
entry in \eqref{eq3.25a}. Hence the quotient $\mathfrak{A}/\cc$ is
a vector bundle of rank $g-2$ over $A$. Since \eqref{eq3.25b} is
equivariant and $\cc$ acts freely on $\mathfrak{A}$, $$\left[
p_*(\cE nd_0\cV\otimes K_X)-p_*\cH om(\cV,K_X)\right]/\cc$$ is a
vector bundle of rank $g$ over $\mathfrak{A}/\cc$.
\end{proof}

Consequently the E-polynomial of the locus of stable pairs of Type
IV is \begin{equation}\label{eq3.27} 2^{2g}\cdot (uv)^{2g-2}\cdot
((uv)^{g-1}-1)\frac{(uv)^g-1}{uv-1}. \end{equation}

\vspace{.5cm}\textbf{3.3. Unstable case}
\bigskip

Suppose $V$ is an unstable rank 2 bundle with trivial determinant.
Then there is a unique line subbundle $L$ of $V$ with maximal
degree and an exact sequence
\begin{equation}\label{eq3.30}
0\to L\to V\to L^{-1}\to 0 \qquad \deg L=d>0
\end{equation}
Our goal in this subsection is to find the subvariety of $\bM^s$
parameterizing stable pairs $(V,\phi)$ with $V$ as in
\eqref{eq3.30} for each $d>0$. If $d>g-1$, then $\deg
(K_XL^{-2})<0$ and $\Hom (L, K_XL^{-1})=0$. This means that $L$ is
preserved by any $\phi\in H^0(\End _0V\otimes K_X)$ and hence
$(V,\phi)$ is never stable. From now on, we let $1\le d\le g-1$.

\begin{proposition} \label{prop3.5} For $1\le d\le g-1$, the
Hodge-Deligne polynomial of the locus of stable pairs $(V,\phi)$
with $V$ as in \eqref{eq3.30} is
\[
(uv)^{3g-3}\cdot E(\tilde{S}^{2g-2-2d}X)
\]
where $\tilde{S}^{2g-2-2d}X$ is a $2^{2g}$-fold covering of the
symmetric product $S^{2g-2-2d}X$ of $X$.\end{proposition}

Let $\Pic^d(X)\to \Pic^{2g-2-2d}(X)$ be the map $L\mapsto
K_XL^{-2}$. This is obviously a $2^{2g}$-fold covering. Let
$P^r_d$ be the locus in $\Pic^{2g-2-2d}(X)$ of line bundles $\xi$
satisfying $h^0(\xi)=r+1$ and $\tilde{P}^r_d$ be the inverse image
of $P^r_d$ in $\Pic^d(X)$. Then $\tilde{P}^r_d$ parameterizes line
bundles $L$ of degree $d$ with $h^1(L^2)=h^0(K_XL^{-2})=r+1$. Let
$S^{2g-2-2d}_rX$ be the inverse image of $P^r_d$ by the
Abel-Jacobi map $S^{2g-2-2d}X\to \Pic^{2g-2-2d}X$ and
$\tilde{S}^{2g-2-2d}_rX$ be the fibre product of $\tilde{P}^r_d$
and $S^{2g-2-2d}_rX$ over $P^r_d$. Then we have a commutative
square
\begin{equation}\label{eq3.31}
\xymatrix{
\tilde{S}^{2g-2-2d}_rX\ar[r]^(.6){\tilde{g}}\ar[d]_{\tilde{f}} &
\tilde{P}^r_d\ar[d]^f\\
S^{2g-2-2d}_rX\ar[r]_(.6)g & P^r_d. }\end{equation} The vertical
maps are $2^{2g}$-fold coverings and the horizontal maps are
$\pp^r$-bundles.

Let $\cL\to \tilde{P}^r_d\times X$ be the restriction of a
universal line bundle over $\Pic^d(X)\times X$. Then there is a
line bundle $\cM\to P^r_d\times X$ such that
\begin{equation}\label{eq3.311}
(f\times 1_X)^*\cM\cong K_X\cL^{-2}\end{equation} by the universal
property of $\Pic(X)$. If we let $\overline{\pi}:P^r_d\times X\to
P^r_d$ be the projection, then
$S^{2g-2-2d}_rX=\pp(\overline{\pi}_*\cM)$ and hence we have an
injection
\begin{equation}\label{eq3.32}
\cO\to g^*\overline{\pi}_*\cM\otimes \cO(1)
\end{equation}
over $S^{2g-2-2d}_rX$. Note that
\begin{equation}\label{eq3.33}
\tilde{f}^*g^*\overline{\pi}_*\cM=\tilde{g}^*f^*\overline{\pi}_*\cM
=\tilde{g}^*\pi_*(f\times 1_X)^*\cM=\tilde{g}^*\pi_*K_X\cL^{-2}
\end{equation}
by \eqref{eq3.311} where $\pi:\tilde{P}^r_d\times X\to
\tilde{P}^r_d$ is the projection. The sheaf $R^1\pi_*\cL^2$ is a
vector bundle of rank $r+1$ and there is a perfect pairing
\begin{equation}\label{eq3.34}
\tilde{g}^*R^1\pi_*\cL^2\otimes \tilde{g}^*\pi_*K_X\cL^{-2}
=\tilde{g}^*(R^1\pi_*\cL^2\otimes \pi_*K_X\cL^{-2})\to
\tilde{g}^*(R^1\pi_*K_X)\cong \cO
\end{equation}
by the Serre duality. Combining \eqref{eq3.32}, \eqref{eq3.33}
and \eqref{eq3.34}, we get a surjective map of vector bundles
\[
\tilde{g}^*R^1\pi_*\cL^2\to \cO(1).
\]
Let $A$ be the kernel of this map which is a vector bundle over
$\tilde{S}^{2g-2-2d}_rX$ of rank $r$. Let $\cL^\#$ be the
pull-back of $\cL$ to $A\times X$ via the composition $A\to
\tilde{S}^{2g-2-2d}_rX\to \tilde{P}^r_d$. It is well-known that
there is a universal extension bundle over $R^1\pi_*\cL^2$ of
$\cL^{-1}$ by $\cL$ and hence we have a rank 2 vector bundle
$\cV^\#$ which fits into an exact sequence
\begin{equation}\label{eq3.35} 0\to \cL^\#\to \cV^\#\to
(\cL^\#)^{-1}\to 0\end{equation} over $A\times X$.

To incorporate the Higgs field, we consider Hitchin's sequence
(\cite[(3.7)]{Hit1}) \begin{equation}\label{eq3.36} 0\to \cH
om(\cV^\#,\cL^\#K_X)\to \cE nd_0\cV^\#\otimes K_X\to
K_X(\cL^\#)^{-2}\to 0
\end{equation}
over $A\times X$. Let $p:A\times X\to A$ be the projection. Then
we have an exact sequence
\begin{equation}\label{eq3.37} 0\to p_*\cH
om(\cV^\#,\cL^\#K_X)\to p_*\cE nd_0\cV^\#\otimes K_X\to
p_*K_X(\cL^\#)^{-2}\to R^1p_*\cH om(\cV^\#,\cL^\#K_X).
\end{equation}
As $\cL^\#$ is the pull-back of $\cL\to \tilde{P}^r_d\times X$,
$p_*K_X(\cL^\#)^{-2}$ is the pull-back of
$\tilde{g}^*\pi_*(K_X\cL^{-2})$. By \eqref{eq3.32} and
\eqref{eq3.33}, there is an injection
\begin{equation}\label{eq3.38}
\cO(-1)\to p_*K_X(\cL^\#)^{-2}.
\end{equation}

From \eqref{eq3.35} we have an exact sequence
\[
0\to K_X(\cL^\#)^{2}\to \cH om(\cV^\#,\cL^\#K_X)\to K_X \to 0
\]
and a long exact sequence
\begin{equation}\label{eq3.40}
0\to p_*K_X(\cL^\#)^{2}\to p_*\cH om(\cV^\#,\cL^\#K_X)\to p_*K_X
\to \end{equation} $$\to R^1p_*K_X(\cL^\#)^{2} \to R^1p_*\cH
om(\cV^\#,\cL^\#K_X)\to R^1p_*K_X\to 0.$$ By the Serre duality,
$R^1p_*K_X(\cL^\#)^{2}=0$ and thus $$R^1p_*\cH
om(\cV^\#,\cL^\#K_X)\cong R^1p_*K_X \cong \cO.$$

By the definition of $A$, the composition of \eqref{eq3.38} with
the last map in \eqref{eq3.37} is zero. Let us consider the
commutative diagram \small
\begin{equation}\label{eq3.39}
\xymatrix{ p_*\cH om(\cV^\#,\cL^\#K_X)\ar[r]& p_*\cE
nd_0\cV^\#\otimes K_X\ar[r] &p_*K_X(\cL^\#)^{-2}\ar[r] & R^1p_*\cH
om(\cV^\#,\cL^\#K_X)\\
p_*\cH om(\cV^\#,\cL^\#K_X)\ar[r]\ar[u]^= & B\ar[r]\ar[u]
&\cO(-1)\ar[u]\ar[r] & 0\ar[u] }\end{equation} \normalsize where
$B$ is the fiber product of $\cO(-1)$ and $p_*\cE
nd_0\cV^\#\otimes K_X$ over $p_*K_X(\cL^\#)^{-2}$. From
\eqref{eq3.40}, $p_*\cH om(\cV^\#,\cL^\#K_X)$ is a vector bundle
of rank $h^0(L^2K_X)+h^0(K_X)=2g+2d-1$ for $L\in \tilde{P}^r_d$
and hence $B$ is a vector bundle over $A$ of rank $2g+2d$.

There is an obvious family $(\cV^\dagger,\Phi^\dagger)$ of Higgs
pairs parameterized by $B$, namely the restriction of the
tautological family of Higgs pairs parameterized by $p_*(\cE
nd_0\cV^\#\otimes K_X)$. A member
$(V,\phi)=(\cV^\dagger,\Phi^\dagger)|_{s}$ for $s\in B$ lying over
$L$ is stable if and only if $\phi$ does not preserve $L$, i.e.
$s$ is not mapped to the zero section of $\cO(-1)$ via the bottom
right horizontal map in \eqref{eq3.39}. Hence the locus of stable
pairs in $B$ is
\[
B-p_*\cH om(\cV^\#,\cL^\#K_X)
\]
which is a $\cc^{2g+2d-1}$-bundle over a $\cc^*\times
\cc^r$-bundle over $\tilde{S}^{2g-2-2d}_rX$. By construction, it
is clear that the bundles are all Zariski locally trivial.

Finally we consider the isomorphism classes of stable pairs in the
family parameterized by $B-p_*\cH om(\cV^\#,\cL^\#K_X)$. Fix $L\in
\tilde{P}^r_d$ and $D\in S^{2g-2-2d}_rX$ such that $K_XL^{-2}\cong
\cO(D)$, i.e. $(L,D)\in \tilde{S}^{2g-2-2d}_rX$. There is an
action of $\cc^*$ on the fiber of $$\left[B-p_*\cH
om(\cV^\#,\cL^\#K_X)\right]\to \tilde{S}^{2g-2-2d}_rX$$ over
$(L,D)$ as follows. Let $(U_i)$ be a sufficiently fine open cover
of $X$. Then for any stable pair $(V,\phi)$ in the fiber, the
transition matrices for $V$ and $\phi$ can be written as
\[
T_{ij}=\left(\begin{matrix} \lambda_{ij} & \rho_{ij}\\
0& \lambda_{ij}^{-1}\end{matrix}\right) \quad \text{ and } \quad
\phi|_{U_i}=\left(\begin{matrix} a_i & b_i\\
c_i& -a_i\end{matrix}\right)
\]
with $c=(c_i)\in H^0(K_XL^{-2})$ whose divisor is
$\mathrm{div}(c)=D$. Then for $t\in\cc^*$, the diagonal matrix
with diagonal entries $(t,t^{-1})$ acts on $T_{ij}$ and
$\phi|_{U_i}$ by conjugation:
\[
t\cdot T_{ij} = \left(\begin{matrix} \lambda_{ij} & t^2\rho_{ij}\\
0& \lambda_{ij}^{-1}\end{matrix}\right),\quad t\cdot \phi|_{U_i}=
\left(\begin{matrix} a_i & t^2b_i\\
t^{-2}c_i& -a_i\end{matrix}\right)
\]
Hence $\cc^*$ acts on the fiber of $\cO(-1)$ and $A$ over
$\tilde{S}^{2g-2-2d}_rX$ with weights $-2$ and $2$ respectively.
But the quotient of $\cc^*\times \cc^r$ by the action of $\cc^*$
with weights $-2$ and $2$ is exactly $\cc^r$. Hence
$$\left[B-p_*\cH
om(\cV^\#,\cL^\#K_X)\right]/\cc^* \cong p_*\cH
om(\cV^\#,\cL^\#K_X)$$ which is a $\cc^{2g+2d-1}$-bundle over $A$.

Next the additive group $(H^0(L^2),+)$ acts on $p_*\cH
om(\cV^\#,\cL^\#K_X)$ as follows. Let $(\mu_i)$ be a cocycle
representing a class in $H^0(L^2)$. Then $\mu_i$ acts on $T_{ij}$
and $\phi|_{U_i}$ by conjugation:
\[
\left(\begin{matrix} 1 & \mu_i\\
0& 1\end{matrix}\right)T_{ij}\left(\begin{matrix} 1 & -\mu_i\\
0& 1\end{matrix}\right)=T_{ij}\]
\begin{equation}\label{eq3.40a}\left(\begin{matrix} 1 & \mu_i\\
0& 1\end{matrix}\right)\left(\begin{matrix} a_i & b_i\\
c_i& -a_i\end{matrix}\right)\left(\begin{matrix} 1 & -\mu_i\\
0& 1\end{matrix}\right)
=\left(\begin{matrix} a_i-c_i\mu_i & b_i+2a_i\mu_i-c\mu_i^2\\
c_i& c_i\mu_i-a_i\end{matrix}\right).\end{equation}

Because $\mathrm{div} (c)=D$, we have an exact sequence
\[\xymatrix{
0\ar[r]& L^2\ar[r]^c &K_X\ar[r] &K_X|_D\ar[r]&0 }\] and hence
$c:H^0(L^2)\to H^0(K_X)$ is injective. Exactly as in the proof of
Lemma \ref{lem3.4} for the Type IV case, the quotient of $p_*K_X$
by the free linear action $a\mapsto a-c\mu$ of $H^0(L^2)$ as in
the first entry of \eqref{eq3.40a} is a vector bundle of rank
$g-h^0(L^2)$ and the second nonzero map in \eqref{eq3.40}
\[
p_*\cH om(\cV^\#,\cL^\#K_X) \to p_*K_X
\]
is equivariant. Hence the quotient of $p_*\cH
om(\cV^\#,\cL^\#K_X)$ by $H^0(L^2)$ is a vector bundle of rank
$h^0(K_XL^2)=g+2d-1$ over a vector bundle of rank
$$g-h^0(L^2)=g-(2d-g+1+r+1)=2g-2d-r-2$$ over $A$ because $h^1(L^2)=r+1$.
Since $A$ is a Zariski locally trivial $\cc^r$-bundle over
$\tilde{S}^{2g-2-2d}_rX$, the E-polynomial of the locus of stable
pairs $(V,\phi)$ with $V$ as in \eqref{eq3.30} and $h^1(L^2)=r+1$,
is
\[
(uv)^{g+2d-1}(uv)^{2g-2d-r-2} (uv)^r E(\tilde{S}^{2g-2-2d}_rX)=
(uv)^{3g-3} E(\tilde{S}^{2g-2-2d}_rX).\] Summing up for all $r$,
we get
\[
(uv)^{3g-3}E(\tilde{S}^{2g-2-2d}X).\] This completes the proof of
Proposition \ref{prop3.5}.

\begin{corollary} The E-polynomial of the locus of stable pairs
with $V$ unstable is
\begin{equation}\label{eq3.50}
(uv)^{3g-3}\cdot\sum_{d=1}^{g-1}E(\tilde{S}^{2g-2-2d}X).
\end{equation}
\end{corollary}

By mimicking Hitchin's computation in \cite[\S7]{Hit1}, we see
that
\[
E(\tilde{S}^{n}X)=E({S}^{n}X)+(2^{2g}-1)\sum_{r+s=n}
(-1)^{r+s}\binom{g-1}{r}\binom{g-1}{s} u^rv^s
\]
\begin{equation}\label{eq3.51}
=
\mathrm{Coeff}_{x^n}\left[\frac{(1-ux)^g(1-vx)^g}{(1-x)(1-uvx)}
+(2^{2g}-1)(1-ux)^{g-1}(1-vx)^{g-1}
\right]
\end{equation}
\[
= \mathrm{Coeff}_{x^n}(1-ux)^{g-1}(1-vx)^{g-1}
\left[\frac{x(1-u)(1-v)}{(1-x)(1-uvx)} +2^{2g} \right].
\]
Observe that
\begin{equation}\label{eq3.52}
\sum_{d=1}^{g-1}\mathrm{Coeff}_{x^{2g-2-2d}}(1-ux)^{g-1}(1-vx)^{g-1}
\end{equation}
$$=\frac12\left[(1-u)^{g-1}(1-v)^{g-1}+(1+u)^{g-1}(1+v)^{g-1}-2(uv)^{g-1}
\right].$$

We keep using Hitchin's method of calculation. We also have
\begin{equation}\label{eq3.53}
\sum_{d=1}^{g-1}\mathrm{Coeff}_{x^{2g-2-2d}}
 \frac{x(1-ux)^{g-1}(1-vx)^{g-1}}{(1-x)(1-uvx)}
\end{equation}
\[
=\sum_{d=1}^{g-1}\mathrm{Res}_{x=0}
 \frac{(1-ux)^{g-1}(1-vx)^{g-1}}{x^{2g-2-2d}(1-x)(1-uvx)}
\]
\[
=\sum_{d=1}^{\infty}\mathrm{Res}_{x=0}
 \frac{(1-ux)^{g-1}(1-vx)^{g-1}}{x^{2g-2-2d}(1-x)(1-uvx)}
\]
\[
=\mathrm{Res}_{x=0}
 \frac{(1-ux)^{g-1}(1-vx)^{g-1}}{x^{2g-4}(1-x^2)(1-x)(1-uvx)}.
\]
As $x\to \infty$, the function is close to $(uv)^{g-2}/x^2$. By
Cauchy's residue theorem,
\[
\mathrm{Res}_{x=0}=-\left(\mathrm{Res}_{x=1}+\mathrm{Res}_{x=-1}
+\mathrm{Res}_{x=(uv)^{-1}}\right).
\]
The residue at the simple pole $x=-1$ is
\[
\frac{(1+u)^{g-1}(1+v)^{g-1}}{4(1+uv)}.
\]
The residue at the simple pole $x=(uv)^{-1}$ is
\[
-\frac{(uv)^{g-1}(1-u)^{g-1}(1-v)^{g-1}}{(uv-1)^2(uv+1)}.
\]
The residue at the double pole $x=1$ is
\[
-\frac{g-1}2\frac{(u+v-2uv)(1-u)^{g-2}(1-v)^{g-2}}{1-uv}\]
\[-\frac{4g-7}4
\frac{(1-u)^{g-1}(1-v)^{g-1}}{1-uv}+\frac{uv(1-u)^{g-1}(1-v)^{g-1}}{2(uv-1)^2}.
\]

Therefore the E-polynomial of the locus of stable pairs $(V,\phi)$
with $V$ unstable is
\begin{equation}\label{eq3.54}
(uv)^{3g-3}\cdot\sum_{d=1}^{g-1}E(\tilde{S}^{2g-2-2d}X)
\end{equation}
$$=2^{2g-1}(uv)^{3g-3}
\left[(1-u)^{g-1}(1-v)^{g-1}+(1+u)^{g-1}(1+v)^{g-1}-2(uv)^{g-1}
\right]$$
$$+(uv)^{3g-3}(1-u)(1-v)\big[
-\frac{(1+u)^{g-1}(1+v)^{g-1}}{4(1+uv)} +
\frac{(uv)^{g-1}(1-u)^{g-1}(1-v)^{g-1}}{(uv-1)^2(uv+1)}$$
\[
+\frac{g-1}2\frac{(u+v-2uv)(1-u)^{g-2}(1-v)^{g-2}}{1-uv}+\frac{4g-7}4
\frac{(1-u)^{g-1}(1-v)^{g-1}}{1-uv}\]
\[-\frac{uv(1-u)^{g-1}(1-v)^{g-1}}{2(uv-1)^2} \big]
\]
from \eqref{eq3.51}, \eqref{eq3.52}, and \eqref{eq3.53}.

\vspace{.5cm}\textbf{3.4. Hodge-Deligne polynomial of $\bM^s$}
\bigskip

So far we considered all possible types of stable pairs $(V,\phi)$
and computed the E-polynomials of the corresponding loci in
$\bM^s$. By adding up \eqref{eq3.1}, \eqref{eq3.10},
\eqref{eq3.17}, \eqref{eq3.23}, \eqref{eq3.27} and \eqref{eq3.54},
we obtain the following.
\begin{theorem}\label{thm1}
\[ E(\bM^s)=
(uv)^{3g-3}
\frac{(1-u^2v)^g(1-uv^2)^g-(uv)^{g+1}(1-u)^g(1-v)^g}{(1-uv)(1-(uv)^2)}\]
$$- (uv)^{3g-3} \frac12(\frac{(1-u)^g(1-v)^g}{1-uv}+\frac{(1+u)^g(1+v)^g}{1+uv})
$$
\[
+(uv)^g\cdot  \big(\frac12 (1-u)^g(1-v)^g+\frac12 (1+u)^g(1+v)^g
\big)\frac{((uv)^g-1)((uv)^{g-1}-1)}{(uv)^2-1}\]
$$+(uv)^{g+1}\cdot \big(\frac12 (1-u)^g(1-v)^g-\frac12 (1+u)^g(1+v)^g
\big) \frac{((uv)^{g-1}-1)((uv)^{g-2}-1)}{(uv)^2-1}
$$
\[
+(uv)^{2g-1} \frac{((uv)^{g-2}-1)((uv)^{g-1}-1)}{uv-1}
[(1-u)^g(1-v)^g-2^{2g}]
\]
$$
+ 2^{2g}\cdot (uv)^{2g-2}\cdot
((uv)^{g-1}-1)\frac{((uv)^{g-1}-1)((uv)^g-1)}{uv-1}.
$$
$$+2^{2g-1}(uv)^{3g-3}
\left[(1-u)^{g-1}(1-v)^{g-1}+(1+u)^{g-1}(1+v)^{g-1}-2(uv)^{g-1}
\right]$$
$$+(uv)^{3g-3}(1-u)(1-v)\big[
-\frac{(1+u)^{g-1}(1+v)^{g-1}}{4(1+uv)} +
\frac{(uv)^{g-1}(1-u)^{g-1}(1-v)^{g-1}}{(uv-1)^2(uv+1)}$$
\[
+\frac{g-1}2\frac{(u+v-2uv)(1-u)^{g-2}(1-v)^{g-2}}{1-uv}+\frac{4g-7}4
\frac{(1-u)^{g-1}(1-v)^{g-1}}{1-uv}\]
\[-\frac{uv(1-u)^{g-1}(1-v)^{g-1}}{2(uv-1)^2} \big].
\]
\end{theorem}

\section{Kirwan's desingularization}

In this section, we show that $\bM$ is desingularized by three
blow-ups by Kirwan's algorithm for desingularizations \cite{k2}.
We will see that the singularities of $\bM$ are identical to those
of the moduli space of rank 2 semistable sheaves with Chern
classes $c_1=0$ and $c_2=2g$ on a K3 surface with generic
polarization as studied in \cite{ogrady}. O'Grady constructed the
Kirwan desingularization by three blow-ups and we use his
arguments.

We first collect some of Simpson's results. Let $N$ be a
sufficiently large integer and $p=2N+2(1-g)$. Then we have the
following.\bigskip

\begin{enumerate}
\item[1.] \cite[Theorem 3.8]{simp}\\ There is a quasi-projective
scheme $Q$ representing the moduli functor which parameterizes the
isomorphism classes of triples $(V,\phi,\alpha)$ where $(V,\phi)$
is a semistable Higgs pair with $\det V\cong \cO_X$, $\mathrm{tr}
\phi=0$ and $\alpha$ is an isomorphism
\[
\alpha:\cc^p\to H^0(X,V\otimes \cO(N)).
\]
\bigskip

\item[2.] \cite[Theorem 4.10]{simp}\\ Fix $x\in X$. Let $\tilde{Q}$
be the frame bundle at $x$ of the universal bundle restricted to
$x$. Then the action of $GL(p)$ lifts to $\tilde{Q}$ and $SL(2)$
acts on the fibers of $\tilde{Q}\to Q$ in an obvious fashion.
Every point of $\tilde{Q}$ is stable with respect to the free
action of $GL(p)$ and
\[\bR=\tilde{Q}/GL(p)\]
represents the moduli functor which parameterizes triples
$(V,\phi,\beta)$ where $(V,\phi)$ is a semistable Higgs pair with
$\det V\cong \cO_X$, $\mathrm{tr} \phi=0$ and $\beta$ is an
isomorphism
\[
\beta:V|_x\to \cc^2.
\]
\bigskip

\item[3.] \cite[Theorem 4.10]{simp}\\
Every point in $\bR$ is semistable with respect to the (residual)
action of $SL(2)$. The closed orbits in $\bR$ correspond to
polystable pairs, i.e. $(V,\phi)$ is stable or \[
(V,\phi)=(L,\psi)\oplus (L^{-1},-\psi)\] for $L\in \Pic^0(X)$ and
$\psi\in H^0(K_X)$. The set $\bR^s$ of stable points with respect
to the action of $SL(2)$  is exactly the locus of stable pairs.
\bigskip

\item[4.] \cite[Theorem 4.10]{simp}\\
The good quotient $\bR\git SL(2)$ is $\bM$.

\bigskip

\item[5.] \cite[Theorem 11.1]{simp}\\
$\bR$ and $\bM$ are both irreducible normal quasi-projective
varieties.

\bigskip

\item[6.] \cite[\S10]{simp}\\
Let $A^i$ (resp. $A^{i,j}$) be the sheaf of smooth $i$-forms
(resp. $(i,j)$-forms) on $X$. For a polystable Higgs pair
$(V,\phi)$, consider the complex
\begin{equation}\label{e4.1}
0\to \End_0V\otimes A^0\to \End_0V\otimes A^1\to \End_0V\otimes
A^2\to 0\end{equation} whose differential is given by
$D''=\overline{\partial}+\phi$. Because $A^1=A^{1,0}\oplus
A^{0,1}$ and $\phi$ is of type $(1,0)$, we have an exact sequence
of complexes with \eqref{e4.1} in the middle
\begin{equation}\label{e4.2}\xymatrix{
& 0\ar[d] & 0\ar[d] &0\ar[d] &    \\
0\ar[r] & 0\ar[r]\ar[d] &\End_0V\otimes
A^{1,0}\ar[r]^{\overline{\partial}}\ar[d] &\End_0\otimes
A^{1,1}\ar[r]\ar[d]^= &0\\
0\ar[r]& \End_0V\otimes A^0\ar[r]^{D''}\ar[d]_= & \End_0V\otimes
A^1\ar[r]^{D''}\ar[d] & \End_0V\otimes A^2\ar[r]\ar[d] & 0\\
0\ar[r] & \End_0V\otimes A^{0,0}\ar[r]^{\overline{\partial}}\ar[d]
&
\End_0V\otimes A^{0,1}\ar[r]\ar[d] & 0\ar[r]\ar[d] & 0 \\
& 0 & 0  &0 }\end{equation} This gives us a long exact sequence
\begin{equation}\label{e4.3}
\xymatrix{ 0\ar[r] & T^0\ar[r] &H^0(\End_0V)\ar[r]^(.42){[\phi,-]}
& H^0(\End_0V\otimes K_X)\ar[r]& }\end{equation}
\[\xymatrix{ \ar[r]& T^1\ar[r] &H^1(\End_0V)\ar[r]^(.42){[\phi,-]}
& H^1(\End_0V\otimes K_X)\ar[r]  &T^2\ar[r] &0 }
\]
where $T^i$ is the $i$-th cohomology of \eqref{e4.1}. The Zariski
tangent space of $\bM$ at polystable $(V,\phi)$ is isomorphic to
$T^1$.
\bigskip

\item[7.] \cite[Theorem 10.4]{simp}\\
Using the above notation, let $C$ be the quadratic cone in $T^1$
defined by the map
\begin{equation}\label{e4.6}
T^1\to T^2
\end{equation} which sends a $\End_0V$-valued 1-form $\eta$ to
$[\eta,\eta]$. Let $y=(V,\phi,\beta)\in \bR$ be a point with
closed orbit. Then the formal completion ${(\bR,y)}^\wedge$ is
isomorphic to the formal completion ${(C\times
\mathfrak{h}^\perp,0)}\wedge$ where $\mathfrak{h}^\perp$ is the
perpendicular space to the image of $T^0\to H^0(\End_0V)\to
sl(2)$. Furthermore, if we let $Y$ be the \'etale slice at $y$ of
the $SL(2)$-orbit in $\bR$, then
\[
(Y,y)^\wedge \cong (C,0)^\wedge.
\]
\bigskip

\item[8.] \cite[Lemma 10.7]{simp}\\
The dimension of the Zariski tangent space to $\bR$ at
$y=(V,\phi,\beta)$ is
\[
\dim T^1+3-\dim T^0
\]
and $\dim T^0=\dim T^2$ by the Serre duality and \eqref{e4.3}.

\end{enumerate}
\bigskip

By Riemann-Roch and \eqref{e4.3}, we have
\begin{equation}\label{e4.5}
\dim T^1= \chi(\End_0\otimes K_X)-\chi(\End_0V)\end{equation}
\[+\dim T^0+\dim T^2=6g-6+2\dim T^0.\]
Therefore, $\bR$ is smooth at $y=(V,\phi,\beta)$ if and only if
$T^0=0$. When $(V,\phi)$ is stable, there is no section of
$H^0(\End_0 V)$ commuting with $\phi$ and so $T^0=0$. Hence the
stable locus $\bR^s$ is smooth and so is the orbit space $\bM^s$.

The complement of $\bR^s$ in $\bR$ consists of 4 subvarieties
parameterizing strictly semistable pairs $(V,\phi)$ of the
following 4 types respectively:
\begin{enumerate}
\item[(i)] $(L,0)\oplus (L,0)$ for $L\cong L^{-1}$
\item[(ii)] nontrivial extension of $(L,0)$ by $(L,0)$ for $L\cong
L^{-1}$
\item[(iii)] $(L,\psi)\oplus (L^{-1},-\psi)$ for $(L,\psi)\ncong (L^{-1},-\psi)$
\item[(iv)] nontrivial extension of $(L^{-1},-\psi)$ by $(L,\psi)$
for $(L,\psi)\ncong (L^{-1},-\psi)$
\end{enumerate}
where $L\in \Pic^0(X)=:J$ and $\psi\in H^0(K_X)$. Higgs pairs of
type (ii) and (iv) are not polystable and their orbits do not
appear in $\bM$.

Let us first consider the locus of type (i). It is obvious from
definition that the loci of type (i) in $\bM$ and in $\bR$ are
both isomorphic to the $\zz_2$-fixed point set $J_0\cong
\zz_2^{2g}$ in the Jacobian $J$ by the involution $L\mapsto
L^{-1}$. From \eqref{e4.3} we have isomorphisms
\[T^0\cong H^0(\End_0V)\cong sl(2)\]
\[T^2\cong H^1(\End_0V\otimes K_X)\cong H^1(K_X)\otimes sl(2)\cong
sl(2)\] and an exact sequence
\begin{equation}\label{e4.7}
0\to H^0(\End_0V\otimes K_X)\to T^1\to H^1(\End_0V)\to 0.
\end{equation}
Choose a splitting of \eqref{e4.7} $$T^1\cong H^0(\End_0V\otimes
K_X)\oplus H^1(\End_0V)\cong H^0(K_X)\otimes sl(2)\oplus
H^1(\cO_X)\otimes sl(2).$$ The quadratic map \eqref{e4.6} is just
the Lie bracket of $sl(2)$ coupled with the perfect pairing
\[ H^0(K_X)\otimes H^1(\cO_X)\to H^1(K_X).\]
It is easy to check that this coincides with the quadratic map
$\overline{\Upsilon}$ in \cite[p.65]{ogrady} and thus the
singularity of $\bR$ along the locus $\zz_2^{2g}$ of type (i) is
the same as the deepest singularity in O'Grady's case
\cite{ogrady}. Moreover, the actions of $SL(2)$ on the quadratic
cones are identical. In particular, the singularity of the locus
$\zz_2^{2g}$ in $\bM$ of Higgs pairs $(L,0)\oplus (L,0)$ with
$L\cong L^{-1}$ is the hyperk\"ahler quotient
$$\mathbb{H}^g\otimes sl(2)/\!/\!/ SL(2)$$
where $\mathbb{H}$ is the division algebra of quaternions.

Next let us look at the locus of type (iii). It is clear that the
locus of type (iii) in $\bM$ is isomorphic to
$$J\times_{\zz_2}H^0(K_X)-\zz_2^{2g}\cong T^*J/\zz_2-\zz_2^{2g}$$
where $\zz_2$ acts on $J$ by $L\mapsto L^{-1}$ and on $H^0(K_X)$
by $\psi\mapsto -\psi$. The locus of type (iii) in $\bR$ is a $\pp
SL(2)/\cc^*$-bundle over $T^*J/\zz_2-\zz_2^{2g}$ and in particular
it is smooth. As in the type (i) case, by using \eqref{e4.3} and
\eqref{e4.6}, it is straightforward to show that the singularity
along the locus of type (iii) in $\bR$ is the same as the
singularity along $\Sigma_Q^0$ in O'Grady's case
\cite[\S1.4]{ogrady}. In particular, the singularity along the
locus in $\bM$ of type (iii) is the hyperk\"ahler quotient
\[ \mathbb{H}^{g-1}\otimes T^*\cc/\!/\!/\cc^*\]
and we have a stratification of $\bM$;
\begin{equation}\label{e4.15}
\bM=\bM^s\sqcup (T^*J/\zz_2-\zz_2^{2g})\sqcup
\zz_2^{2g}.\end{equation}

Since we have identical singularities and stratification as in
O'Grady's case in \cite{ogrady}, we can copy his arguments almost
line by line to construct the Kirwan desingularization $\hat{\bM}$
of $\bM$ and study its blow-downs. So we skip the details but give
only a brief outline.

We blow up $\bR$ first along the locus of type (i) and then along
the locus of type (iii). Then we delete the unstable part with
respect to the action of $SL(2)$. After these two blow-ups the
loci of types (ii) and (iv) become unstable by \cite{k2} (or more
explicitly by \cite[Lemma 1.7.4]{ogrady}) and thus they are
deleted anyway.  Let $\bS^{ss}$ (resp. $\bS^s$) be the open subset
of semistable (resp. stable) points after the two blow-ups. Then
we have \begin{enumerate}
\item[(a)] $\bS^{ss}=\bS^s$,
\item[(b)] $\bS^s$ is smooth.\end{enumerate}
In particular, $\bS^s/ SL(2)$ has at worst orbifold singularities.
When $g=2$, this is already smooth. When $g\ge 3$, by blowing up
$\bS^s$ one more time along the locus of points with stabilizers
larger than the center $\zz_2$ of $SL(2)$, we obtain a variety
$\hat{\bS}$ whose orbit space
$$\hat{\bM}:=\hat{\bS}/SL(2)$$
is a smooth variety obtained by blowing up $\bM$ first along
$\zz_2^{2g}$, second along the proper transform of $T^*J/\zz_2$
and third along a nonsingular subvariety lying in the proper
transform of the exceptional divisor of the first blow-up. Let
$$\pi:\hat{\bM}\to \bM$$
be the composition of the three blow-ups. We call $\hat{\bM}$ the
Kirwan desingularization of $\bM$. Along the line of O'Grady's
\cite{ogrady}, we obtain the following.

\begin{theorem}
\begin{enumerate}
\item For $g=2$, $\hat{\bM}$ can be blown down to give us a
symplectic desingularization of $\bM$.
\item For $g\ge 3$, $\hat{\bM}$ can be blown down twice to give us
another smooth model of $\bM$, which we call the O'Grady
desingularization.
\item The three exceptional divisors  $D_1,D_2,D_3$ of
$\pi:\hat{\bM}\to \bM$ coming from the three blow-ups are smooth
and normal crossing. The discrepancy divisor is
\begin{equation}\label{e4.16}
K_{\hat{\bM}}=K_{\hat{\bM}}-\pi^* K _{\bM}
=(6g-7)D_1+(2g-4)D_2+(4g-6)D_3.
\end{equation}
\end{enumerate}\end{theorem}
Note that $K_{\bM}=0$ since $\bM^s$ is hyperk\"ahler \cite{Hit1}.

It is possible to extract explicit descriptions of the divisors
$D_1,D_2,D_3$ from \cite{ogrady} as follows. (See also
\cite[Proposition 3.2]{CK1}, \cite[Proposition 3.6]{CK2}.) Let
$\hat{\pp}^5$ be the blow-up of $\pp^5$ (projectivization of the
space of $3\times 3$ symmetric matrices) along $\pp^2$ (the locus
of rank 1 matrices). Let $(\cc^{2g},\omega)$ be a symplectic
vector space and let $\Gr^{\omega}(k,2g)$ be the Grassmannian of
$k$-dimensional subspaces of $\cc^{2g}$, isotropic with respect to
the symplectic form $\omega$ (i.e. the restriction of $\omega$ to
the subspace is zero). Let $I_{2g-3}$ denote the incidence variety
given by
\[ I_{2g-3}=\{(p,H)\in \pp^{2g-3}\times \breve{\pp}^{2g-3}| p\in
H\}.
\]
Then we have the following.

\begin{proposition}\label{pro4.1}
Let $g\geq 3$. \begin{enumerate} \item[(1)] $D_1$ is the disjoint
union of $2^{2g}$ copies of a $\hat{\pp}^5$-bundle over
$\Gr^\omega(3,2g)$.

\item[(2)] $D_2^0$ is a free $\zz_2$-quotient of a
$I_{2g-3}$-bundle over $ T^*J-J_0$.

\item[(3)] $D_3$ is the disjoint
union of $2^{2g}$ copies of a $\pp^{2g-4}$-bundle over a
$\pp^2$-bundle over $\Gr^\omega(2,2g)$.

\item[(4)] $D_{12}$ is the disjoint
union of $2^{2g}$ copies of a $\pp^2$-bundle over a $\pp^2$-bundle
over $\Gr^\omega(3,2g)$.

\item[(5)] $D_{23}$ is the disjoint
union of $2^{2g}$ copies of a $\pp^{2g-4}$-bundle over a  $
\pp^1$-bundle over $\Gr^\omega(2,2g)$.

\item[(6)] $D_{13}$ is the disjoint
union of $2^{2g}$ copies of a $ \pp^2$-bundle over a
$\pp^2$-bundle over $\Gr^\omega(3,2g)$.

\item[(7)] $D_{123}$ is the disjoint
union of $2^{2g}$ copies of a $\pp^1$-bundle over a $\pp^2$-bundle
over $\Gr^\omega(3,2g)$. \end{enumerate} All the bundles above are
Zariski locally trivial.
\end{proposition}

\section{Stringy E-function of $\bM$}

In this section we compute the stringy E-function of $\bM$ by
using Proposition \ref{pro4.1} and \eqref{e4.16} and show that
there does not exist a symplectic desingularization of $\bM$ for
$g\ge 3$.

By Simpson's theorem, $\bM$ is an irreducible normal variety with
Gorenstein singularities because $K_{\bM}=0$. By \eqref{e4.16},
the singularities are canonical and the stringy E-function of
$\bM$ is a well-defined rational function of $u,v$. From
(\ref{eqn:stringy E-function}) and \eqref{e4.16}, the stringy
E-function $E_{st}(\bM)=E_{st}(\bM;u,v)$ of $\bM$ is
\begin{equation}\label{eqn:stringy E-function of M_c}
E(\bM^s;u,v)+E(D^0_1;u,v){\textstyle\frac{1-uv}{1-(uv)^{6g-6}}}
+E(D^0_2;u,v){\textstyle\frac{1-uv}{1-(uv)^{2g-3}}}\end{equation}
\[
+E(D^0_3;u,v) {\textstyle\frac{1-uv}{1-(uv)^{4g-5}}}
+E(D^0_{12};u,v){\textstyle\frac{1-uv}{1-(uv)^{6g-6}}\frac{1-uv}{1-(uv)^{2g-3}}}
\]
\[+E(D^0_{23};u,v){\textstyle\frac{1-uv}{1-(uv)^{2g-3}}\frac{1-uv}{1-(uv)^{4g-5}}}
+E(D^0_{13};u,v){\textstyle\frac{1-uv}{1-(uv)^{4g-5}}\frac{1-uv}{1-(uv)^{6g-6}}}
\]
\[ +E(D^0_{123};u,v){\textstyle\frac{1-uv}{1-(uv)^{6g-6}}
\frac{1-uv}{1-(uv)^{2g-3}}\frac{1-uv}{1-(uv)^{4g-5}}} .\]

From Proposition \ref{pro4.1} and the identity (\cite[Lemma
3.1]{CK1}, \cite[Lemma 4.1]{CK2})
\[ E(\Gr^\omega(k,2g);u,v)=\prod_{1\leq i\leq k} \frac{1-(uv)^{2g-2k+2i}}{1-(uv)^i}, \]
we obtain the following. (See \cite[Corollary 3.3]{CK1} or
\cite[Corollary 4.2]{CK2}.)

\begin{proposition}\label{eqn:computation of stringy E-function}
$$ E(D_1;u,v) = 2^{2g}\cdot \Bigl({\textstyle
\frac{1-(uv)^6}{1-uv}-\!\frac{1-(uv)^3}{1-uv}+\!\bigl(\frac{1-(uv)^3}{1-uv}\bigr)^2}\Bigr)
\cdot \prod_{1\leq i\leq 3}\! \Bigl({\textstyle
\frac{1-(uv)^{2g-6+2i}}{1-(uv)^i}}\Bigr),$$

$$E(D_3;u,v)   = 2^{2g}\cdot {\textstyle
\frac{1-(uv)^{2g-3}}{1-uv}\cdot\frac{1-(uv)^3}{1-uv}} \cdot
\prod_{1\leq i\leq 2}\Bigl({\textstyle \frac{1-(uv)^{2g-4+2i}}
{1-(uv)^i}}\Bigr), $$

$$ E(D_{12};u,v)   = 2^{2g}\cdot \Bigl({\textstyle
\frac{1-(uv)^3}{1-uv}}\Bigr)^2\cdot \prod_{1\leq i\leq
3}\Bigl({\textstyle \frac{1-(uv)^{2g-6+2i}}{1-(uv)^i}}\Bigr), $$

$$ E(D_{23};u,v)  = 2^{2g}\cdot {\textstyle
\frac{1-(uv)^{2g-3}}{1-uv}\cdot\frac{1-(uv)^2}{1-uv}} \cdot
\prod_{1\leq i\leq 2}\Bigl({\textstyle
\frac{1-(uv)^{2g-4+2i}}{1-(uv)^i}}\Bigr), $$

$$E(D_{13};u,v)   = 2^{2g}\cdot {\textstyle \frac{
1-(uv)^3}{1-uv}\cdot\frac{1-(uv)^{2g-4}}{1-uv}} \cdot \prod_{1\leq
i\leq 2}\Bigl({\textstyle
\frac{1-(uv)^{2g-4+2i}}{1-(uv)^i}}\Bigr), $$
$$ E(D_{123};u,v)   =2^{2g}\cdot{\textstyle
\frac{1-(uv)^2}{1-uv}\cdot\frac{1-(uv)^{2g-4}}{1-uv}}\cdot
\prod_{1\leq i\leq 2}\Bigl({\textstyle
\frac{1-(uv)^{2g-4+2i}}{1-(uv)^i}}\Bigr).$$
\end{proposition}

For $D_2^0$, observe that $H^*_c(D_2^0)$ is the $\zz_2$-invariant
part of
$$H^*(I_{2g-3})\otimes H^*_c(T^*J-J_0)$$ and hence we have
$$E(D_2^0;u,v)=E(I_{2g-3};u,v)^+\cdot \left(E(T^*J;u,v)^+-2^{2g}\right) +
E(I_{2g-3};u,v)^-\cdot E(T^*J;u,v)^-$$ where $E(Z;u,v)^+$ (resp.
$E(Z;u,v)^-$) denotes the E-polynomial of the $\zz_2$-invariant
(resp. anti-invariant) part of $H^*_c(Z)$. By elementary
computation (\cite[\S5]{CK1} or \cite[Lemma 4.3]{CK2}), we have
\[
E(I_{2g-3};u,v)^+=\frac{(1-(uv)^{2g-2})(1-(uv)^{2g-3})}{(1-uv)(1-(uv)^2)}
\]
\[
E(I_{2g-3};u,v)^-=uv\cdot
\frac{(1-(uv)^{2g-2})(1-(uv)^{2g-3})}{(1-uv)(1-(uv)^2)}.
\]
It is also elementary that
\[
E(T^*J;u,v)^+=\frac12 (uv)^g\cdot \left[
(1-u)^g(1-v)^g+(1+u)^g(1+v)^g\right]
\]
\[
E(T^*J;u,v)^-=\frac12 (uv)^g\cdot \left[
(1-u)^g(1-v)^g-(1+u)^g(1+v)^g\right].
\]
Therefore we have
\begin{equation}\label{e5.2}
E(D_2^0;u,v)=(uv)^g\cdot
\frac{(1-(uv)^{2g-2})(1-(uv)^{2g-3})}{(1-uv)(1-(uv)^2)}\end{equation}
\[
\times \left[ \frac12 (1+uv) (1-u)^g(1-v)^g+ \frac12 (1-uv)
(1+u)^g(1+v)^g -2^{2g}\right].
\]

By the additive property of the Hodge-Deligne polynomials we have
\[
E(D_1^0;u,v)= 2^{2g}\cdot ((uv)^5-(uv)^2)\cdot \prod_{1\leq i\leq
3} \frac{1-(uv)^{2g-6+2i}}{1-(uv)^i}
\]
\[
E(D_3^0;u,v)= 2^{2g}\cdot (uv)^{2g-2}\cdot \prod_{1\leq i\leq 2}
\frac{1-(uv)^{2g-4+2i}}{1-(uv)^i}
\]
\[
E(D_{12}^0;u,v)= 2^{2g}\cdot (uv)^2\cdot (1+uv+(uv)^2)\cdot
\prod_{1\leq i\leq 3} \frac{1-(uv)^{2g-6+2i}}{1-(uv)^i}
\]
\[
E(D_{23}^0;u,v)= 2^{2g}\cdot (uv)^{2g-4}\cdot (1+uv)\cdot
\prod_{1\leq i\leq 2} \frac{1-(uv)^{2g-4+2i}}{1-(uv)^i}
\]
\[
E(D_{13}^0;u,v)= 2^{2g}\cdot (uv)^2\cdot (1+uv+(uv)^2)\cdot
\prod_{1\leq i\leq 3} \frac{1-(uv)^{2g-6+2i}}{1-(uv)^i}
\]
\[
E(D_{123}^0;u,v)= 2^{2g}\cdot (1+uv)\cdot (1+uv+(uv)^2)\cdot
\prod_{1\leq i\leq 3} \frac{1-(uv)^{2g-6+2i}}{1-(uv)^i}.
\]

By direct computation with \eqref{eqn:stringy E-function of M_c},
\eqref{e5.2} and the above, we obtain that $E_{st}(\bM)-E(\bM^s)$
is equal to
\begin{equation}\label{eqlastminus}
(uv)^g\cdot \frac{1-(uv)^{2g-2}}{1-(uv)^2} \cdot \left[ \frac12
(1+uv) (1-u)^g(1-v)^g+ \frac12 (1-uv) (1+u)^g(1+v)^g
-2^{2g}\right]
\end{equation}
\[
+ 2^{2g}\cdot
\frac{(1-(uv)^{2g-2})(1-(uv)^{2g})}{1-(uv)^{4g-5}}\cdot \big[
\frac{ 1-(uv)^{8g-10}}{(1-(uv)^{2g-3})(1-(uv)^{6g-6})}
\]
\[
+ \frac{(uv)^2(1-(uv)^{2g-4})
(1-(uv)^{6g-8})}{(1-(uv)^{2})(1-(uv)^{2g-3})(1-(uv)^{6g-6})}
+\frac{(uv)^{2g-2}}{1-(uv)^{2}} \big].
\]

By Theorem \ref{thm1} and the above, we proved the following.
\begin{theorem}
\[
E_{st}(\bM)= (uv)^{3g-3}
\frac{(1-u^2v)^g(1-uv^2)^g-(uv)^{g+1}(1-u)^g(1-v)^g}{(1-uv)(1-(uv)^2)}\]
$$- (uv)^{3g-3} \frac12(\frac{(1-u)^g(1-v)^g}{1-uv}+\frac{(1+u)^g(1+v)^g}{1+uv})
$$
\[
+(uv)^g\cdot  \big(\frac12 (1-u)^g(1-v)^g+\frac12 (1+u)^g(1+v)^g
\big)\frac{((uv)^g-1)((uv)^{g-1}-1)}{(uv)^2-1}\]
$$+(uv)^{g+1}\cdot \big(\frac12 (1-u)^g(1-v)^g-\frac12 (1+u)^g(1+v)^g
\big) \frac{((uv)^{g-1}-1)((uv)^{g-2}-1)}{(uv)^2-1}
$$
\[
+(uv)^{2g-1} \frac{((uv)^{g-2}-1)((uv)^{g-1}-1)}{uv-1}
[(1-u)^g(1-v)^g-2^{2g}]
\]
$$
+ 2^{2g}\cdot (uv)^{2g-2}\cdot
((uv)^{g-1}-1)\frac{((uv)^{g-1}-1)((uv)^g-1)}{uv-1}.
$$
$$+2^{2g-1}(uv)^{3g-3}
\left[(1-u)^{g-1}(1-v)^{g-1}+(1+u)^{g-1}(1+v)^{g-1}-2(uv)^{g-1}
\right]$$
$$+(uv)^{3g-3}(1-u)(1-v)\big[
-\frac{(1+u)^{g-1}(1+v)^{g-1}}{4(1+uv)} +
\frac{(uv)^{g-1}(1-u)^{g-1}(1-v)^{g-1}}{(uv-1)^2(uv+1)}$$
\[
+\frac{g-1}2\frac{(u+v-2uv)(1-u)^{g-2}(1-v)^{g-2}}{1-uv}+\frac{4g-7}4
\frac{(1-u)^{g-1}(1-v)^{g-1}}{1-uv}\]
\[-\frac{uv(1-u)^{g-1}(1-v)^{g-1}}{2(uv-1)^2} \big]
\]
\[
+(uv)^g\cdot \frac{1-(uv)^{2g-2}}{1-(uv)^2} \cdot \left[ \frac12
(1+uv) (1-u)^g(1-v)^g+ \frac12 (1-uv) (1+u)^g(1+v)^g
-2^{2g}\right]
\]
\[
+ 2^{2g}\cdot
\frac{(1-(uv)^{2g-2})(1-(uv)^{2g})}{1-(uv)^{4g-5}}\cdot \big[
\frac{ 1-(uv)^{8g-10}}{(1-(uv)^{2g-3})(1-(uv)^{6g-6})}
\]
\[
+ \frac{(uv)^2(1-(uv)^{2g-4})
(1-(uv)^{6g-8})}{(1-(uv)^{2})(1-(uv)^{2g-3})(1-(uv)^{6g-6})}
+\frac{(uv)^{2g-2}}{1-(uv)^{2}} \big].
\]

\end{theorem}

By taking the limit $u,v\to 1$, we obtain the stringy Euler number
$$e_{st}(\bM)=\lim_{u,v\to 1}E_{st}(\bM;u,v)$$
as follows.
\begin{corollary}
\[ e_{st}(\bM)= 2^{2g}\frac{3g-3}{2g-3}.\]
\end{corollary}

In particular, the stringy Euler number is never an integer for
$g\ge 4$. When $g=3$, one can check directly that $E_{st}(\bM)$
(or \eqref{eqlastminus}) is not a polynomial.
\begin{corollary}
There does not exist a symplectic desingularization of $\bM$ for
$g\ge 3$. \end{corollary}
\begin{proof} Since $\bM^s$ is hyperk\"ahler, the canonical bundle
$K_\bM$ is trivial. If there were a symplectic desingularization
$\mathfrak{M}$ of $\bM$, then $K_{\mathfrak{M}}=0$ by definition
and hence the resolution would be crepant. But in that case,
$E_{st}(\bM)$ has to be equal to the Hodge-Deligne polynomial of
$\mathfrak{M}$ which is a polynomial with integer coefficients.
This contradicts the fact that $e_{st}(\bM)$ is not an integer for
$g\ge 4$ and that $E_{st}(\bM)$ is not a polynomial for
$g=3$.\end{proof}

\bibliographystyle{amsplain}

\end{document}